# A new perspective in linear Cauchy Elasticity: variational minimum principles for statics, dynamics, and heterogeneous materials

**Amit Acharya**

Department of Civil & Environmental Engineering, and Center for Nonlinear Analysis, Carnegie Mellon University, Pittsburgh, USA

A variational minimum principle for linear elastodynamics of a possibly heterogeneous material without a stored energy function is developed. It involves a change of variables to dual fields, and results in a degenerate elliptic Euler-Lagrange system, even when the primal elastodynamics is hyperbolic. Uniqueness assertions for the dual dynamic and static problems and implications of the degenerate ellipticity are sketched. Some implications pertaining to heterogeneous materials and ones with indefinite elastic moduli are discussed.

## 1 Introduction

This work describes an application to Linear Cauchy Elasticity of a technique (Acharya 2023a; Acharya 2023b; Acharya 2025) for associating a family of variational minimum (i.e., convex) principles to a system of equations with boundary and initial conditions when appropriate. The given 'primal' equations may not have a variational structure in terms of the primal fields or variables they are originally posed in, e.g., elasticity without a strain energy function. The procedure involves an adapted change of variables to a new set of 'dual' fields in terms of which a functional can be posed whose Euler-Lagrange equations become the primal set of equations, interpreted through the change of variables mapping involved. This approach has been further developed and demonstrated through computation in (Kouskiya and Acharya 2024; Kouskiya and Acharya 2025; Acharya and Sengupta 2024; Singh et al. 2024; Kouskiya et al. 2025; Sukumar and Acharya 2025), with rigorous results developed in (Singh et al. 2024; Acharya, Ginster, and Sengupta 2024; Acharya and Ginster 2025; Kouskiya et al. 2025; Vorotnikov and Acharya 2025; Acharya, Stroffolini, et al. 2024). Here, we work out the details of the scheme in the context of the governing system being that of linear elastodynamics of a generally heterogeneous material possibly without a strain energy function. As is well understood, accounting for initial conditions of an initial value problem is problematic in a variational principle—including Hamilton's celebrated principle—see, e.g., (Gurtin 1964) for a discussion; Hamilton's principle also has a problem with causality due to the need for an acausal final-time boundary condition on position, cf. (Willis 1981; Galley 2013; Horowitz and Rothkopf 2026). The way out of this difficulty in linear problems has been to resort to working in the Laplace transform domain (Gurtin 1964; Willis 1981), or at 'fixed frequency' (Milton et al. 2009), i.e., for elastodynamic fields satisfying a special ansatz for its time-dependent part, assumed separable from its space-dependent part; also see (Horowitz and Rothkopf 2026) in the context of Hamiltonian particle mechanics using an adaptation of the Schwinger-Keldysh formalism which involves a doubling of variables and a backward-in-time path. The approach adopted here works in real space-time domains without involving convolutions, assumptions of separability, or a doubling of variables in which the original system is posed, and works as well for nonlinear and dissipative time-dependent equations as shown in (Kouskiya and Acharya 2025; Sukumar and Acharya 2025; Vorotnikov and Acharya 2025; Singh et al. 2024).

The work (Singh et al. 2024) develops the dual approach for nonlinear elasticity, focusing primarily on elastostatics; the form of the dual functional is developed for quadratic nonlinearity in stress. The emphasis in this work is to explicitly work out the case of linear elasticity in full detail, which is also of interest to a much broader class of researchers interested in its theoretical and/or practical aspects. As examples, the specialization to linear elastodynamics greatly facilitates the explorations of the type provided in Sections 3 and 4. For comparison to other variational techniques for linear elastodynamics in the literature, to our knowledge this work provides the first convex variational principle for linear elastodynamics (with a strain energy function and without) in the time domain without approximation, recovering field equations, boundary conditions, *and* initial conditions without any loss of causality. In situations when a stored energy function exists, the recovery of initial conditions, a causal formulation, and the convexity of the variational principle for the whole range of positive-definite to indefinite stored energy function are some of the main contributions of this work.

An outline of the paper is as follows: in Section 2 the variational principle for elastodynamics is developed. In Section 3 a uniqueness result for the Euler-Lagrange system of the variational principle is studied; some discussion of its ellipticity properties (in space-time domains) is also provided. In Section 4 implications of the scheme for heterogeneous/composite materials are discussed. Some concluding remarks are recorded in Section 5.

## 2 A convex variational principle for linear elastodynamics

Given a 3-dimensional domain $\Omega$ with boundary $\partial\Omega = \overline{\partial\Omega_s \cup \partial\Omega_u}$ the closure of a union of disjoint sets (with $\partial\Omega_u$ containing at least three non collinear points), and an interval of time $[0,T]$, consider the system of linear elastodynamics whose given elastic moduli and density fields, $(x,t) \mapsto C(x), \rho(x)$ can vary in space. $C$ has minor symmetries and may or may not have major symmetries; we will use the notation $(C^\top)_{ijkl} = C_{klij}$. The fields $u, v$ will denote the displacement and velocity, and $E$ and $W$ the symmetric and antisymmetric parts of the displacement gradient. When using direct notation a " $\cdot$ " will denote the inner product of two vectors and a " : " that of two second-order tensors. For fourth-order tensors, $A, B, C$, $(AB)_{ijkl} = A_{ijrm}B_{rmkl}$ and $(ABC)_{ijkl} = A_{ijrm}B_{rmps}C_{pskl}$. For this work, it will be convenient to work with a non-dimensional representation of the equations of linear elasticity. Using $E^*$ as a physical scale for stress[1], $L^*$ as some physical scale for length[2], and $T^*$ as a scale for time[3], the equations of linear elastodynamics can be written as

$$\partial_{\tilde{t}}\tilde{u}_i - \tilde{v}_i = 0 \qquad \text{on } \Omega \times (0,T) \tag{1a}$$

$$\tilde{\rho}\partial_{\tilde{t}}\tilde{v}_i - \partial_{\tilde{x}_j}(\tilde{C}_{ijkl}E_{kl}) - \tilde{b}_i = 0 \qquad \text{on } \Omega \times (0,T) \tag{1b}$$

$$\partial_{\tilde{x}_l}\tilde{u}_k - E_{kl} - W_{kl} = 0 \qquad \text{on } \Omega \times (0,T) \tag{1c}$$

$$\tilde{u}_i = \tilde{u}_i^{(b)} \qquad \text{on } \partial\Omega_u \times (0,T) \tag{1d}$$

$$(\tilde{C}_{ijkl}E_{kl})n_j = \tilde{t}_i \qquad \text{on } \partial\Omega_s \times (0,T) \tag{1e}$$

$$\tilde{u}_i = \tilde{u}_i^{(0)} \qquad \text{on } \Omega \times \{0\} \tag{1f}$$

$$\tilde{v}_i = \tilde{v}_i^{(0)} \qquad \text{on } \Omega \times \{0\}. \tag{1g}$$

where $\tilde{\rho} = \rho L^{*2}/(E^* T^{*2})$, $\tilde{x}_i = x_i/L^*$, $\tilde{t} = t/T^*$, $\tilde{u}_i = u_i/L^*$, $\tilde{v}_i = v_i T^*/L^*$, $\tilde{C}_{ijkl} = C_{ijkl}/E^*$, $\tilde{u}_i^{(0)} = u_i^{(0)}/L^*$, $\tilde{v}_i^{(0)} = v_i^{(0)}T^*/L^*$, $\tilde{u}_i^{(b)} = u_i^{(b)}/L^*$, $\tilde{t}_i = t_i/E^*$, $\tilde{b}_i = b_i L^*/E^*$. Here, $b$ the body force density, $u^{(b)}, t$ the displacement and traction boundary conditions, and $u^{(0)}, v^{(0)}$ the displacement and velocity initial conditions. Henceforth, we work with the system (1) *dropping all tildes* for convenience:

$$
\begin{aligned}
&\partial_t u_i - v_i = 0 && \text{on } \Omega \times (0,T) && \text{(2a)}\\
&\rho \partial_t v_i - \partial_j (C_{ijkl} E_{kl}) - b_i = 0 && \text{on } \Omega \times (0,T) && \text{(2b)}\\
&\partial_l u_k - E_{kl} - W_{kl} = 0 && \text{on } \Omega \times (0,T) && \text{(2c)}\\
&u_i = u_i^{(b)} && \text{on } \partial\Omega_u \times (0,T) && \text{(2d)}\\
&(C_{ijkl} E_{kl}) n_j = t_i && \text{on } \partial\Omega_s \times (0,T) && \text{(2e)}\\
&u_i = u_i^{(0)} && \text{on } \Omega \times \{0\} && \text{(2f)}\\
&v_i = v_i^{(0)} && \text{on } \Omega \times \{0\}. && \text{(2g)}
\end{aligned}
$$

[1] E.g., the Young's modulus for the material.

[2] E.g., dimension of a body, a feature-size like grain size, a wave length of imposed disturbance, or a ratio of an elastic modulus, say $E^*$, and the magnitude of the body force density at some point or its average over a physically relevant region in an infinite homogeneous medium.

[3] E.g., the time for an elastic wave to traverse the length $L^*$ given by $L^*/\sqrt{E^*/\rho^*}$ where $\rho^*$ is a constant mass density scale; or a reciprocal frequency of harmonic loading.

Using the notation $U := (u, v, E, W)$, and $\bar{U} := (\bar{u}, \bar{v}, \bar{E}, \bar{W})$ where $(x,t) \mapsto \bar{U}(x,t)$ are arbitrarily specified functions referred to as *base states* (in space-time domains), we introduce an *auxiliary potential*

$$
\begin{aligned}
H(U, \bar{U}) = \frac{1}{2}\big( &(E_{ij} - \bar{E}_{ij}) A_{ijkl} (E_{kl} - \bar{E}_{kl}) + (W_{kl} - \bar{W}_{kl})(W_{kl} - \bar{W}_{kl}) \\
&+ (u_i - \bar{u}_i)(u_i - \bar{u}_i) + \mathrm{r} (v_i - \bar{v}_i)(v_i - \bar{v}_i) \big)
\end{aligned} \tag{3}
$$

for the system of linear elastodynamics, where $A, \mathrm{r}$ are nondimensional fields, with freedom to be designed as needed. Due to the quadratic form through which it enters, $A$ is assumed to have major symmetry without loss of generality.

Consider now a *pre-dual* functional

$$
\begin{aligned}
\widehat{S}[U, \gamma, \lambda, \mu] = &\int_0^T \int_\Omega \big( -u_i \partial_t \gamma_i - \gamma_i v_i - v_i \partial_t (\rho \lambda_i) + E_{kl} C_{ijkl} \partial_j \lambda_i \\
&\qquad - u_k \partial_l \mu_{kl} - \mu_{kl} (E_{kl} + W_{kl}) - H(U, \bar{U}) \big) \, \mathrm{d}x \, \mathrm{d}t \\
&- \int_0^T \int_\Omega \lambda_i b_i \, \mathrm{d}x \, \mathrm{d}t - \int_\Omega \big( \gamma_i(x,0) u_i^{(0)}(x) + \rho(x) \lambda_i(x,0) v_i^{(0)}(x) \big) \, \mathrm{d}x \\
&- \int_0^T \int_{\partial\Omega_s} \lambda_i t_i \, \mathrm{d}a \, \mathrm{d}t + \int_0^T \int_{\partial\Omega_u} \mu_{kl} u_k^{(b)} n_l \, \mathrm{d}a \, \mathrm{d}t
\end{aligned} \tag{4}
$$

with

$$
\begin{aligned}
&\gamma_i = 0 && \text{on } \Omega \times \{T\} && \text{(5a)}\\
&\lambda_i = 0 && \text{on } \Omega \times \{T\}; && \text{(5b)}\\
&\lambda_i = 0 && \text{on } \partial\Omega_u \times (0,T) && \text{(5c)}\\
&\mu_{kl} n_l = 0 && \text{on } \partial\Omega_s \times (0,T). && \text{(5d)}
\end{aligned}
$$

The fields $D := (\gamma, \lambda, \mu)$ are referred to as *dual* fields. The values (right hand-sides) of the dual (Dirichlet) space-time boundary conditions (5) can be arbitrarily specified, without loss of generality (note that in such cases no extra terms are, nevertheless, included in the pre-dual functional which is free to design.)

We refer to the first space-time bulk integrand in Equation (4) as the Lagrangian

$$
\mathcal{L}(U, \mathcal{D}, \bar{U}) \qquad \text{where} \qquad \mathcal{D} := (D, \nabla D, \partial_t D). \tag{6}
$$

Defining the DtP *(dual-to-primal) map* $(\mathcal{D}, \bar{U}) \mapsto \hat{U}(\mathcal{D}, \bar{U})$ as the solution of

$$
\frac{\partial \mathcal{L}}{\partial U}(U, \mathcal{D}, \bar{U}) = 0
$$

for $U$ in terms of $(\mathcal{D}, \bar{U})$ so that

$$
\frac{\partial \mathcal{L}}{\partial U}(\hat{U}(\mathcal{D}, \bar{U}), \mathcal{D}, \bar{U}) = 0, \tag{7}
$$

we have

$$\hat{u}_m - \bar{u}_m = (-\partial_t \gamma_m - \partial_l \mu_{ml}) \tag{8a}$$
$$\hat{v}_m - \bar{v}_m = \frac{1}{\mathrm{r}}(-\gamma_m - \rho \partial_t \lambda_m) \tag{8b}$$
$$\hat{E}_{rs} - \bar{E}_{rs} = A^{-1}_{rsmn}(C_{ijmn}\partial_j \lambda_i - \mu_{(mn)}) \tag{8c}$$
$$\hat{W}_{mn} - \bar{W}_{mn} = -\mu_{[mn]}, \tag{8d}$$

where square brackets around a pair of indices represents antisymmetrization e.g., ($A_{[ij]} = \frac{1}{2}(A_{ij} - A_{ji})$) and ordinary brackets, symmetrization ($A_{(ij)} = \frac{1}{2}(A_{ij} + A_{ji})$). Also, for any second-order tensor, say $D$, $D^{(s)}$, $D^{(a)}$ will represent its symmetric and antisymmetric parts, respectively.

Now define the *dual* functional, $S$, by substituting $\hat{U}$ for $U$ in the pre-dual functional $\widehat{S}$:

$$\begin{aligned} S[D] = & \int_0^T \int_\Omega \mathcal{L}(\hat{U}(\mathcal{D}, \bar{U}), \mathcal{D}, \bar{U})\, \mathrm{d}x\, \mathrm{d}t \\ & - \int_0^T \int_\Omega \lambda_i b_i \, \mathrm{d}x\, \mathrm{d}t - \int_\Omega (\gamma_i(x,0) u_i^{(0)}(x) + \rho(x)\lambda_i(x,0) v_i^{(0)}(x))\, \mathrm{d}x \\ & - \int_0^T \int_{\partial\Omega_s} \lambda_i t_i \, \mathrm{d}a\, \mathrm{d}t + \int_0^T \int_{\partial\Omega_u} \mu_{kl} u_k^{(b)} n_l \, \mathrm{d}a\, \mathrm{d}t \end{aligned} \tag{9}$$

with the essential space-time boundary conditions (5). Then, noting Equation (7) and that the Lagrangian (6) is affine in $\mathcal{D}$ by construction, it is easily seen that the Euler-Lagrange equations of the functional $S$ is the governing set of equations and side conditions of linear Cauchy elastodynamics (2) with the replacement $U \longrightarrow \hat{U}$, *regardless of the choices of the specified functions*

$$(x,t) \quad \mapsto \quad A(x,t), \quad \mathrm{r}(x), \quad \bar{u}(x,t), \quad \bar{v}(x,t), \quad \bar{E}(x,t), \quad \bar{W}(x,t),$$

*subject to* $\mathrm{r} \neq 0$ *and* $A$ *merely invertible on the space of symmetric second-orders tensors.*

Let us now assume that $\mathrm{r} > 0$ and $A$ is positive-definite on the space of symmetric second order tensors and note that

$$\widehat{S}[U,D] = \int_0^T \int_\Omega \mathcal{L}(U, \mathcal{D}, \bar{U})\, \mathrm{d}x\, \mathrm{d}t + \text{ boundary terms involving dual and specified fields.}$$

Then, due to the lack of any differential constraints on $U$ in the definition of $\mathcal{L}$,

$$\int_0^T \int_\Omega \mathcal{L}(\hat{U}(\mathcal{D}, \bar{U}), \mathcal{D}, \bar{U})\, \mathrm{d}x\, \mathrm{d}t = \int_0^T \int_\Omega \sup_U \mathcal{L}(U, \mathcal{D}, \bar{U})\, \mathrm{d}x\, \mathrm{d}t = \sup_U \int_0^T \int_\Omega \mathcal{L}(U, \mathcal{D}, \bar{U})\, \mathrm{d}x\, \mathrm{d}t.$$

Therefore,

$$S[D] = \sup_U \widehat{S}[U,D],$$

and because $\widehat{S}$ is affine in $D$ by construction, and therefore convex in $D$, defining our problem statement as finding the infimum of $S$, i.e.,

$$\inf_D S[D] = \inf_D \sup_U \widehat{S}[U,D]; \qquad D^* = \operatorname*{argmin}_D S[D]$$

corresponds to a minimization problem for a *convex* functional, solving formally the problem of linear Cauchy Elastodynamics (statics is a special case) through the use of the DtP map, including for heterogeneous materials.

To obtain the explicit form of the dual functional, the (dual) Lagrangian may be expressed, using the DtP map (8), as

$$
\begin{aligned}
&\mathcal{L}(\hat{U},\mathcal{D},\bar{U}) + H(\hat{U},\bar{U}) = \\
&\quad = (\hat{u}_i - \bar{u}_i)\hat{u}_i + \mathsf{r}(\hat{v}_i - \bar{v}_i)\hat{v}_i + \hat{E}_{mn}A_{mnrs}(\hat{E}_{rs} - \bar{E}_{rs}) + (\hat{W}_{mn} - \bar{W}_{mn})\hat{W}_{mn} \\
&\quad = 2H(\hat{U},\bar{U}) + (\hat{u}_i - \bar{u}_i)\bar{u}_i + \mathsf{r}(\hat{v}_i - \bar{v}_i)\bar{v}_i + A_{mnrs}(\hat{E}_{rs} - \bar{E}_{rs})\bar{E}_{mn} + (\hat{W}_{mn} - \bar{W}_{mn})\bar{W}_{mn} \\
&\quad = (\partial_t\gamma_m + \partial_l\mu_{ml})(\partial_t\gamma_m + \partial_l\mu_{ml}) + \tfrac{1}{\mathsf{r}}(\gamma_m + \rho\partial_t\lambda_m)(\gamma_m + \rho\partial_t\lambda_m) \\
&\qquad + (C_{rsmn}\partial_s\lambda_r - \mu_{(mn)})A^{-1}_{mnij}A_{ijkl}A^{-1}_{klpq}(C_{rspq}\partial_s\lambda_r - \mu_{(pq)}) + \mu_{[mn]}\mu_{[mn]} \\
&\qquad\quad - (\partial_t\gamma_m + \partial_l\mu_{ml})\bar{u}_m - (\gamma_m + \rho\partial_t\lambda_m)\bar{v}_m \\
&\qquad\qquad + (C_{ijmn}\partial_j\lambda_i - \mu_{(mn)})\bar{E}_{mn} - \mu_{[mn]}\bar{W}_{mn}.
\end{aligned} \tag{10}
$$

Thus,

$$
\begin{aligned}
S[D] = \frac{1}{2}\int_0^T\!\!\int_\Omega &\big(|\partial_t\gamma + \operatorname{div}\mu|^2 + \tfrac{1}{\mathsf{r}}|\gamma + \rho\partial_t\lambda|^2 \\
&\quad + (C^\top\nabla\lambda - \mu^{(s)}) : A^{-1}(C^\top\nabla\lambda - \mu^{(s)}) + \mu^{(a)} : \mu^{(a)}\big)\,\mathrm{d}x\,\mathrm{d}t \\
&- \int_0^T\!\!\int_\Omega \big((\partial_t\gamma + \operatorname{div}\mu)\cdot\bar{u} + (\gamma + \rho\partial_t\lambda)\cdot\bar{v} + (\mu - C^\top\nabla\lambda) : \bar{E} + \mu^{(a)} : \bar{W}\big)\,\mathrm{d}x\,\mathrm{d}t \\
&- \int_0^T\!\!\int_\Omega \lambda\cdot b\,\mathrm{d}x\,\mathrm{d}t - \int_\Omega \big(\gamma(x,0)\cdot u^{(0)} + \rho\lambda(x,0)\cdot v^{(0)}\big)\,\mathrm{d}x \\
&- \int_0^T\!\!\int_{\partial\Omega_s} \lambda\cdot t\,\mathrm{d}a\,\mathrm{d}t + \int_0^T\!\!\int_{\partial\Omega_u} (\mu n)\cdot u^{(b)}\,\mathrm{d}a\,\mathrm{d}t.
\end{aligned} \tag{11}
$$

with its fields satisfying the following essential boundary conditions on the space-time domain $\Omega \times (0,T)$:

$$\gamma = 0 \quad \text{on } \Omega \times \{T\} \tag{12a}$$
$$\lambda = 0 \quad \text{on } \Omega \times \{T\}; \tag{12b}$$
$$\lambda = 0 \quad \text{on } \partial\Omega_u \times (0,T) \tag{12c}$$
$$\mu n = 0 \quad \text{on } \partial\Omega_s \times (0,T). \tag{12d}$$

For $A$ positive-definite and $\mathsf{r} > 0$, the semi-definiteness of the second variation of the representation (11) of the dual functional furnishes another (direct) proof of its convexity, regardless of the symmetry or positivity of $C$ or $\rho$.

The DtP map for the problem is given by

$$\hat{u} = -(\partial_t\gamma + \operatorname{div}\mu) + \bar{u} \tag{13a}$$
$$\hat{v} = -\tfrac{1}{\mathsf{r}}(\gamma + \rho\partial_t\lambda) + \bar{v} \tag{13b}$$
$$\hat{E} = -A^{-1}(\mu^{(s)} - C^\top\nabla\lambda) + \bar{E} \tag{13c}$$
$$\hat{W} = -\mu^{(a)} + \bar{W}. \tag{13d}$$

The system of linear elastodynamics along with its initial and boundary conditions, expressed in terms of the dual fields, is

$$\partial_t(\bar{u} - (\partial_t\gamma + \operatorname{div}\mu)) - (\bar{v} - \tfrac{1}{\mathsf{r}}(\gamma + \rho\partial_t\lambda)) = 0 \quad \text{on } \Omega \times (0,T) \tag{14a}$$
$$\rho\partial_t(\bar{v} - \tfrac{1}{\mathsf{r}}(\gamma + \rho\partial_t\lambda)) - \operatorname{div}(C(\bar{E} - A^{-1}(\mu^{(s)} - C^\top\nabla\lambda))) - b = 0 \quad \text{on } \Omega \times (0,T) \tag{14b}$$
$$\nabla(\bar{u} - (\partial_t\gamma + \operatorname{div}\mu)) - (\bar{E} - A^{-1}(\mu^{(s)} - C^\top\nabla\lambda)) - (\bar{W} - \mu^{(a)}) = 0 \quad \text{on } \Omega \times (0,T) \tag{14c}$$
$$\bar{u} - (\partial_t\gamma + \operatorname{div}\mu) = u^{(b)} \quad \text{on } \partial\Omega_u \times (0,T) \tag{14d}$$
$$(C(\bar{E} - A^{-1}(\mu^{(s)} - C^\top : \nabla\lambda)))n = t \quad \text{on } \partial\Omega_s \times (0,T) \tag{14e}$$
$$\bar{u} - (\partial_t\gamma + \operatorname{div}\mu) = u^{(0)} \quad \text{on } \Omega \times \{0\} \tag{14f}$$
$$\bar{v} - \tfrac{1}{\mathsf{r}}(\gamma + \rho\partial_t\lambda) = v^{(0)} \quad \text{on } \Omega \times \{0\} \tag{14g}$$

with (these space-time b.c.s can be arbitrarily assigned):

$$\gamma = 0 \qquad \text{on } \Omega \times \{T\} \tag{15a}$$
$$\lambda = 0 \qquad \text{on } \Omega \times \{T\} \tag{15b}$$
$$\lambda = 0 \qquad \text{on } \partial\Omega_u \times (0, T) \tag{15c}$$
$$\mu n = 0 \qquad \text{on } \partial\Omega_s \times (0, T). \tag{15d}$$

The following observations are in order:

1. When the base state, $\bar{U}$, is a solution to the system (2), $D = 0$ is a solution to the dual system (14) and an (absolute) minimizer of the dual functional, $S$ in Equation (11). This is an important property of the dual scheme for both theoretical and practical purposes (Vorotnikov and Acharya 2025), as it motivates the reasoning that if a base state $\bar{U}$ is 'close' to a solution of Equation (2), then it may be possible to obtain a solution of Equation (2) by using $D = 0$ as a good guess to the corresponding dual problem (14).
2. The first variation of the dual functional $S$ in Equation (9) imposes the linear elastodynamics equations (2) in the weak form - thus, all jump conditions of the primal problem are preserved by any extremal of the dual variational problem.
3. When solutions of the primal system (2) are unique, *any* solution to the dual system (14), generated from, e.g., different choices of boundary conditions (5) or $(A, \mathrm{r}, \bar{U})$, must generate the unique primal solution through the use of the DtP map. Examples of this fact are provided in (Kouskiya and Acharya 2024) in the context of the heat equation. However, if the primal system is known to have non-unique solutions, the auxiliary potential, $H$ parametrized by $(A, \mathrm{r}, \bar{U})$, can be used to obtain such solutions through the dual scheme in a 'stable' manner—e.g., even when $C$ is indefinite. In other words, the use of the auxiliary potential, $H$, acts as a selection criterion for solutions of the primal problem.
4. Given that the primal problem (2) is an initial value-problem a natural question is whether prescribing final-time Dirichlet boundary conditions on the dual fields interferes in any way with causality of the primal solutions obtained from the above dual approach. Such interference does not arise because the DtP map (8) involves time derivatives of the dual fields, and prescribing the values of dual fields at final-time $T$ allows these time derivatives (at time $T$) to adjust so as to recover the correct primal solution. This feature has been explained and demonstrated in (Acharya 2023b, Section 7) and (Kouskiya and Acharya 2024; Kouskiya and Acharya 2025; Sukumar and Acharya 2025), with an exact solution for the *first-order* initial value problem $\dot{u} = au, u(0) = u_0$ worked out by the dual scheme in (Sukumar and Acharya 2025, Section 3.1).
5. Due to the availability of the functional $S$, a gradient descent algorithm in a fake time-like variable, can be formulated to obtain solutions to the primal system (2), see (Vorotnikov and Acharya 2025, Section 3). The fact that $S$ is convex also guarantees that the norm of its gradient is non-increasing along a gradient flow—this is of great importance in a solution scheme that employs more than one $S$ functional (based on choices of base states) while keeping the gradient continuous at such switches, with the final goal of reaching a vanishing norm of the gradient (a weak solution to the primal equations is then obtained). Such ideas can be useful for iterative schemes for linear elasticity of heterogeneous materials.
6. The point of view adopted here is that all physics is contained in the primal system and the variational scheme is simply a mathematical device to solve that system. Thus, the dual fields are unlikely to have much physical significance—e.g., they admit acausal boundary conditions. Also, in nonlinear problems the use of a whole sequence of dual functionals, parametrized by base states, is employed to obtain a single primal solution; physically important quantities are not usually found in such abundance.
7. From a practical standpoint, having to solve a boundary-value-problem on a large space-time domain as might be the case to probe long-time behavior can be prohibitive. Solution of the dual system (14) can be divided into space-time domain slices, with an arbitrary, finite, sub-division of the time interval $[0, T]$. In each such space-time slice the dual problem can be solved, the primal solution recovered at the final time of the slice and used to define the initial conditions for the next slice. This strategy for solving the dual problem has been demonstrated in the solution of Euler's ODE system for motion of a rigid body about a fixed point (with and without damping), the heat equation, and the linear transport equation in (Kouskiya and Acharya 2024), and in the context of the (in)viscid Burgers equation in (Kouskiya and Acharya 2025).

## 3 Formal uniqueness for the dual system (14) and its degenerate ellipticity

For practical purposes related to the use of the dual system (14), it is reassuring to have an assertion of uniqueness of solutions for any specific set of dual Dirichlet boundary conditions (5) and base states $\bar{U}$, at least when the elastic modulus $C$ is positive definite and has major symmetry. We provide such a sufficient condition in this Section. In the static case, uniqueness holds simply with positive definiteness without an assumption of major symmetry and we sketch that proof as well.

Consider two solutions of system (14) for $(\gamma, \lambda, \mu)$ and let their difference be denoted as $(\gamma^*, \lambda^*, \mu^*)$. Then

$$\partial_t(\partial_t\gamma^* + \operatorname{div}\mu^*) - \tfrac{1}{\mathsf{r}}(\gamma^* + \rho\partial_t\lambda^*) = 0 \qquad \text{on } \Omega \times (0,T) \tag{16a}$$

$$\rho\partial_t(\tfrac{1}{\mathsf{r}}(\gamma^* + \rho\partial_t\lambda^*)) - \operatorname{div}(C(A^{-1}(\mu^{*(s)} - C^\top\nabla\lambda^*))) = 0 \qquad \text{on } \Omega \times (0,T) \tag{16b}$$

$$\nabla(\partial_t\gamma^* + \operatorname{div}\mu^*) - A^{-1}(\mu^{*(s)} - C^\top\nabla\lambda^*) - \mu^{*(a)} = 0 \qquad \text{on } \Omega \times (0,T) \tag{16c}$$

Also, the difference fields satisfy the following space-time boundary conditions:

$$\partial_t\gamma^* + \operatorname{div}\mu^* = 0 \qquad \text{on } \partial\Omega_u \times (0,T) \tag{17a}$$

$$(C(A^{-1}(\mu^{*(s)} - C^\top\nabla\lambda^*)))n = 0 \qquad \text{on } \partial\Omega_s \times (0,T) \tag{17b}$$

$$(\partial_t\gamma^* + \operatorname{div}\mu^*) = 0 \qquad \text{on } \Omega \times \{0\} \tag{17c}$$

$$\tfrac{1}{\mathsf{r}}(\gamma^* + \rho\partial_t\lambda^*) = 0 \qquad \text{on } \Omega \times \{0\} \tag{17d}$$

$$\gamma^* = 0 \qquad \text{on } \Omega \times \{T\} \tag{17e}$$

$$\lambda^* = 0 \qquad \text{on } \Omega \times \{T\} \tag{17f}$$

$$\lambda^* = 0 \qquad \text{on } \partial\Omega_u \times (0,T) \tag{17g}$$

$$\mu^* n = 0 \qquad \text{on } \partial\Omega_s \times (0,T). \tag{17h}$$

On forming products of the three sets of equations in (16) with $(\gamma^*, \lambda^*, \mu^*)$, respectively, and using the space-time boundary conditions (17a)-(17d), one obtains,

$$0 = \int_0^T\!\!\int_\Omega \Big(|\partial_t\gamma^* + \operatorname{div}\mu^*|^2 + \frac{1}{\mathsf{r}}|\gamma^* + \rho\partial_t\lambda^*|^2 + (C^\top\nabla\lambda^* - \mu^{*(s)}) : A^{-1}(C^\top\nabla\lambda^* - \mu^{*(s)}) + \mu^{*(a)} : \mu^{*(a)}\Big)\,dx\,dt,$$

and due to the positive-definiteness of $A$ on the space of symmetric tensors and $\mathsf{r} > 0$,

$$\partial_t\gamma^* + \operatorname{div}\mu^* = 0 \tag{18a}$$

$$\gamma^* + \rho\partial_t\lambda^* = 0 \tag{18b}$$

$$C^\top\nabla\lambda^* - \mu^{*(s)} = 0 \tag{18c}$$

$$\mu^{*(a)} = 0 \tag{18d}$$

hold *a.e.* on $\Omega \times (0,T)$. Then Equation (18) (all four equations together) gives

$$\rho\partial_{tt}\lambda^* = \operatorname{div}(C^\top\nabla\lambda^*), \tag{19}$$

which further gives, after forming a scalar product with $\partial_t\lambda^*$ and integrating over $\Omega$, *when $C$ has major symmetry,*

$$\partial_t\int_\Omega\Big(\frac{1}{2}\rho|\partial_t\lambda^*|^2 + \frac{1}{2}\nabla\lambda^* : C^\top\nabla\lambda^*\Big)\,dx - \int_{\partial\Omega}\partial_t\lambda^*\cdot((C^\top\nabla\lambda^*)n)\,da = 0. \tag{20}$$

Now, from Equation (17g) one has that

$$\partial_t \lambda^* = 0 \text{ on } \partial\Omega_u \times (0, T). \tag{21}$$

Assuming that (using Equation (18d) as motivation)

$$\mu^{*(a)} n = 0 \text{ on } \partial\Omega_s \times (0, T) \tag{22}$$

holds, and combining with Equation (17h) gives

$$\mu^{*(s)} n = 0 \text{ on } \partial\Omega_s \times (0, T). \tag{23}$$

*Assuming* as well that (using Equation (18c) as motivation),

$$(C^\top \nabla\lambda^* - \mu^{*(s)}) n = 0 \text{ on } \partial\Omega_s \times (0, T), \tag{24}$$

Equations (21), (23) and (24) imply that the boundary term in Equation (20) vanishes so that

$$\int_{\Omega\times\{T\}} (\rho|\partial_t \lambda^*|^2 + \nabla\lambda^* : C^\top \nabla\lambda^*)\, \mathrm{d}x - \int_{\Omega\times\{t\}} (\rho|\partial_t \lambda^*|^2 + \nabla\lambda^* : C^\top \nabla\lambda^*)\, \mathrm{d}x = 0 \tag{25}$$

for all $t \in (0, T)$. *Assuming* (using Equation (18b) as motivation)

$$\gamma^* + \rho\partial_t \lambda^* = 0 \text{ on } \Omega \times \{T\} \tag{26}$$

and noting Equations (17e) and (17f), the first integral in Equation (25) vanishes and *when* $\rho > 0$ *and* $C$ *is positive-definite* (major symmetry has been assumed), we have that

$$\partial_t \lambda^* = 0 \qquad \text{on } \Omega \times (0, T) \tag{27a}$$

$$\nabla\lambda^{*(s)} = 0 \qquad \text{on } \Omega \times (0, T). \tag{27b}$$

This implies, from Equations (18c) and (18d) (and an assumption of continuous extension to the boundary) that $\mu^* = 0$ on $\bar{\Omega} \times [0, T]$, which combined with Equations (17e) and (18a) gives $\gamma^* = 0$ on $\bar{\Omega} \times [0, T]$. Of course, Equation (27a) along with Equation (17f) implies that $\lambda^* = 0$ on $\bar{\Omega} \times [0, T]$. The (formal) argument above for dual elastodynamics may be summarized as follows:

**Proposition 1** Consider a class of functions $(\gamma, \lambda, \mu)$, $\gamma, \lambda : \bar{\Omega} \times [0, T] \to \mathbb{R}^3$, $\mu : \bar{\Omega} \times [0, T] \to \mathbb{R}^{3\times3}$, such that the difference $(\Delta\gamma, \Delta\lambda, \Delta\mu)$ of any two functions in the class satisfy the conditions (17), (22), (24), and (26) with the replacement $(\gamma^*, \lambda^*, \mu^*) \to (\Delta\gamma, \Delta\lambda, \Delta\mu)$. Furthermore, assume that $\Delta\mu$ is continuous on $\bar{\Omega} \times [0, T]$. When $\rho > 0$ and $C$ is positive-definite with major symmetry, there can be at most one solution in this class of functions to the system (14).

**Uniqueness for dual Elastostatics** In the problem of dual elastostatics, the $\gamma, \gamma^*$ dual fields are redundant, and the difference solution $(\lambda^*, \mu^*)$ satisfies

$$-\operatorname{div}(C(A^{-1}(\mu^{*(s)} - C^\top \nabla\lambda^*))) = 0 \qquad \text{on } \Omega \tag{28a}$$

$$\nabla \operatorname{div} \mu^* - A^{-1}(\mu^{*(s)} - C^\top \nabla\lambda^*) - \mu^{*(a)} = 0 \qquad \text{on } \Omega \tag{28b}$$

with the boundary conditions

$$(C(A^{-1}(\mu^{*(s)} - C^\top \nabla\lambda^*)))n = 0 \qquad \text{on } \partial\Omega_s \tag{29a}$$

$$\lambda^* = 0 \qquad \text{on } \partial\Omega_u \tag{29b}$$

$$\mu^* n = 0 \qquad \text{on } \partial\Omega_s. \tag{29c}$$

Then, following the previous arguments for the dynamic case, one has

$$\operatorname{div} \mu^* = 0 \qquad \text{on } \Omega \tag{30a}$$

$$C^\top \nabla\lambda^* - \mu^{*(s)} = 0 \qquad \text{on } \Omega \tag{30b}$$

$$\mu^{*(a)} = 0 \qquad \text{on } \Omega, \tag{30c}$$

which further implies

$$0 = \operatorname{div}(C^\top \nabla \lambda^*) \qquad \text{on } \Omega. \tag{31}$$

Forming the scalar product with $\lambda^*$, integrating by parts and using Equation (29b) gives

$$-\int_\Omega \nabla \lambda^* : C^\top \nabla \lambda^* \, dx + \int_{\partial\Omega_s} \lambda^* \cdot (C^\top \nabla \lambda^*) n \, da = 0 \tag{32}$$

and *assuming* Equations (30b) and (30c) also hold on $\partial\Omega_s$, and combining with Equation (29c), the boundary term in Equation (32) vanishes.

Then, *assuming C to be positive-definite* (but not necessarily symmetric, a difference from what was used in the proof of uniqueness for the dynamic case), one has $\nabla \lambda^{*(s)} = 0$ on $\Omega$ and combining with Equation (29b), uniqueness for $\lambda$ holds. But then, from Equations (30b) and (30c), uniqueness for $\mu$ holds as well.

The (formal) argument above for dual elastostatics may be summarized as follows:

**Proposition 2** Consider a class of functions $(\lambda, \mu)$, $\lambda : \bar{\Omega} \to \mathbb{R}^3$, $\mu : \bar{\Omega} \to \mathbb{R}^{3\times 3}$, such that the difference $(\Delta\lambda, \Delta\mu)$ of any two functions in the class satisfy the conditions i) (30b) and (30c) on $\partial\Omega_s$ and ii) (29), with the replacement $(\lambda^*, \mu^*) \to (\Delta\lambda, \Delta\mu)$.

When $C$ is positive-definite (major symmetry is not required), there can be at most one solution in this class of functions to the system

$$-\operatorname{div}(C(\bar{E} - A^{-1}(\mu^{(s)} - C^\top \nabla\lambda))) - b = 0 \qquad \text{on } \Omega \tag{33a}$$
$$\nabla(\bar{u} - \operatorname{div}\mu) - (\bar{E} - A^{-1}(\mu^{(s)} - C^\top \nabla\lambda)) - (\bar{W} - \mu^{(a)}) = 0 \qquad \text{on } \Omega \tag{33b}$$
$$\bar{u} - \operatorname{div}\mu = u^{(b)} \qquad \text{on } \partial\Omega_u \tag{33c}$$
$$(C(\bar{E} - A^{-1}(\mu^{(s)} - C^\top : \nabla\lambda)))n = t \qquad \text{on } \partial\Omega_s \tag{33d}$$
$$\lambda = 0 \qquad \text{on } \partial\Omega_u \tag{33e}$$
$$\mu n = 0 \qquad \text{on } \partial\Omega_s \tag{33f}$$

where the last two boundary conditions can be arbitrarily assigned.

Returning to the full dual elastodynamic problem, it is to be noted that while uniqueness in the absence of strict convexity might seem surprising, degenerate elliptic equations are known to have uniqueness in many circumstances, see, e.g. (Punzo and Tesei 2009) (here we deal with a system of equations).

Purely from the standpoint of the operator containing the highest order derivatives of the dual system (14) given by

$$-\begin{pmatrix} \partial_{tt} & 0 & \partial_t \partial_m \\ 0 & \frac{\rho_o^2}{\mathsf{r}}\partial_{tt} & \partial_j (C : A^{-1} : C^\top)_{ijmn} \partial_n \\ \partial_j \partial_t & 0 & \partial_j \partial_m \end{pmatrix} \begin{pmatrix} \gamma_i \\ \lambda_m \\ \mu_{im} \end{pmatrix} \tag{34}$$

with the corresponding quadratic form

$$\frac{1}{2} \int_0^T \int_\Omega (|\partial_t \gamma + \operatorname{div}\mu|^2 + \frac{\rho^2}{\mathsf{r}} |\partial_t \lambda|^2 + \nabla\lambda : (C : A^{-1} : C^\top) \nabla\lambda) \, dx \, dt, \tag{35}$$

the system (14) is at worst degenerate-elliptic since the quadratic form (35) is positive semi-definite, *regardless of whether C is positive-definite* and *ρ is positive* (assuming $A$ is positive definite and $\mathsf{r} > 0$, which are free to choose). This property is true of the full quadratic form in Equation (11) as well.

Thus, it is worthy of note that the degenerate ellipticity of the dual problem holds even when the primal elastodynamic problem is hyperbolic or ill-posed hyperbolic with no continuous dependence on initial data when $C$ is indefinite. Initial demonstrations of success in solving such problems in the context of the linear transport equation, the (in)viscid Burgers equation, and elastodynamics of a bar with non-convex strain energy have been provided in (Kouskiya and Acharya 2024; Kouskiya and Acharya 2025; Singh et al. 2024).

We end this Section by exploring the following question: even strictly convex linear elastodynamics admits stress-waves induced by boundary/initial conditions (or body forces), which travel at the linear elastic wave speed of the material. This naturally must involve a propagating strain discontinuity as well. How can a 'close-to' elliptic PDE formulation recover such intrinsically 'hyperbolic' behavior?

Consider the dual 'action' (11) (for $C_{ijkl} = C_{klij}$ and $C$ positive-definite) and one asks about solutions that have vanishing bulk 'dual energy/action.' This necessarily requires the system (18) to hold (where the $^*$ fields now represents the 0-dual-action solution under consideration), which further implies that Equation (19) holds. But this directly implies, because of its hyperbolic nature, that $\nabla\lambda$ admits rank-one discontinuities in the usual way across the characteristic surfaces of the physical linear elastodynamic equation and, through the DtP map (8), in the physical strain field $\hat{E}$ produced by the dual scheme, *regardless of the chosen* $A$, $\mathrm{r}$, $\bar{U}$. On the other hand, in the static situation there can be no discontinuities in $\lambda$ (for $C$ positive-definite) by the usual properties of elliptic equations. Thus, we conclude that the degenerate ellipticity of the dual system (14) is crucial for recovering physically mandated behavior as per classical linear elasticity, and ellipticity of the dual system for elastodynamics must fail. These are all manifestations of a convex minimum principle which is not strictly convex, see (Acharya 2025, Section 3).

The above arguments lead to a natural conjecture: that the dual formulation, due to its degenerate ellipticity selects only 'good' solutions, i.e. ones with fairly mild discontinuities, for elastodynamics with indefinite $C$ tensor.

# 4 Heterogeneous materials

For application to heterogeneous materials, i.e., $C, \rho$ vary in space, it is useful to write the system (14) in the form (we focus only on the governing equations)

$$-\partial_{tt}\gamma - \partial_t \operatorname{div}\mu + \gamma = -\partial_t\bar{u} + \bar{v} - \tfrac{1}{\mathrm{r}}\gamma - \tfrac{\rho}{\mathrm{r}}\partial_t\lambda + \gamma \tag{36a}$$

$$-\tfrac{\rho^2}{\mathrm{r}}\partial_{tt}\lambda - \operatorname{div}((CA^{-1}C^\top)\nabla\lambda) = -\rho\partial_t\bar{v} + \operatorname{div}(C\bar{E}) + \tfrac{\rho}{\mathrm{r}}\partial_t\gamma - \operatorname{div}((CA^{-1})\mu^{(s)}) + b \tag{36b}$$

$$-\nabla\partial_t\gamma - \nabla\operatorname{div}\mu + \mu = -\nabla\bar{u} + \bar{E} + \bar{W} + \mu^{(s)} - A^{-1}\mu^{(s)} + (A^{-1}C^\top)\nabla\lambda. \tag{36c}$$

In the above, in the first equation we have added a term $\gamma$ to both sides of the equation; for the third equation we have done the same using $\mu^{(s)}$.

With the goal of converting the left-hand-side operator of Equation (36) to a constant coefficient one (the 'comparison medium'), let $S^{(0)}$ be a positive-definite, spatially constant, compliance tensor with major and minor symmetries which is otherwise arbitrary and $\rho_0 > 0$ a spatially constant density field. Define

$$C^{(0)} = (S^{(0)})^{-1}$$

and assume $C$ (and $C^\top$) are invertible on the space of symmetric second order tensors. One then makes the choice

$$A := C^\top S^{(0)} C \qquad \text{and} \qquad \mathrm{r} := \rho^2/\rho_0, \tag{37}$$

resulting in Equation (36) taking the form

$$-\partial_{tt}\gamma - \partial_t \operatorname{div}\mu + \gamma = -\partial_t\bar{u} + \bar{v} - \tfrac{\rho_0}{\rho^2}\gamma - \tfrac{\rho_0}{\rho}\partial_t\lambda + \gamma \tag{38a}$$

$$-\rho_0\partial_{tt}\lambda - \operatorname{div}(C^{(0)}\nabla\lambda) = -\rho\partial_t\bar{v} + \operatorname{div}(C\bar{E}) + \tfrac{\rho_0}{\rho}\partial_t\gamma - \operatorname{div}((C^{(0)}(C^\top)^{-1})\mu^{(s)}) + b \tag{38b}$$

$$-\nabla\partial_t\gamma - \nabla\operatorname{div}\mu + \mu = -\nabla\bar{u} + \bar{E} + \bar{W} + \mu^{(s)} - (C^{-1}C^{(0)}(C^\top)^{-1})\mu^{(s)} + (C^{-1}C^{(0)})\nabla\lambda. \tag{38c}$$

It is worthy of note that obtaining a constant-coefficient highest order operator on the l.h.s. of Equation (38) does not come at the expense of introducing a highest order operator on the

r.h.s. of the equation, as is effectively the case in the conventional method used to solve the linear elastic problem for heterogeneous materials (Hashin and Shtrikman 1962; Willis 1981; Moulinec and Suquet 1994) (we will refer to this as the 'established formulation'). There, the addition arises from the (negative) divergence of the stress polarization tensor in the static case:

$$\operatorname{div}(C\nabla u) = \operatorname{div}(C^{(0)}\nabla u) + \operatorname{div}\tau; \qquad \tau = (C - C^{(0)})\nabla u. \tag{39}$$

To elaborate a bit further, consider the elastostatic case and ignore the base states ($\bar{U}$) for the moment in Equation (38). On defining the polarization tensor here to be

$$\tau = (C^{(0)}(C^{\top})^{-1})\mu^{(s)}, \tag{40}$$

Equation (38b) appears to have a very similar structure to the established formulation for solving the problem for a given stress polarization field, as described in the above references. Here, the polarization depends on the tensor field $\mu$ which, very roughly speaking, from Equation (38c) would seem to have one less spatial derivative than $\nabla\lambda$, the latter being the analogous situation in the case in the established formulation. Thus, it is conceivable that solving for $\nabla\lambda$ and $\mu$ here, using a natural adaptation of the Moulinec-Suquet (Moulinec and Suquet 1994) iterative scheme for solving such integral equations, and through these fields the solution for the physical strain field through the DtP map (8c), has better properties. A similar argument would also be applicable to the problem of heterogeneous elastodynamics.

We note the following features, potentially important for possible practical application:

1. A gradient flow algorithm (Vorotnikov and Acharya 2025) based on the convex dual functional $S$ (11) can be utilized as an alternate methodology to the Moulinec-Suquet integral equation based iterative scheme (Moulinec and Suquet 1994) for obtaining (approximate) solutions to system (2) for a heterogeneous body, employing the choices (37).
2. In either approach, changes of base states in an iterative scheme (Vorotnikov and Acharya 2025, Section 3), based on progressive iterates for the primal fields, can be invoked to enhance convergence to a solution.

## 5 Concluding remarks

A scheme for generating variational minimum principles for Cauchy elastodynamics has been described as an application of a broader program to ease the solution/approximation of difficult (linear and nonlinear) PDE problems (Acharya 2025). This is achieved through a conversion of the question to a *convex* variational one which allows the definition of a weaker notion of solutions (*Variational dual solutions* (Singh et al. 2024; Vorotnikov and Acharya 2025; Acharya, Stroffolini, et al. 2024)) to PDE than weak solutions, thus allowing for the use of tools from the Calculus of Variations and Convex Analysis in PDE. Such solutions have the consistency that when they can be proven to be sufficiently regular, they define genuine weak solutions of the primal PDE problem. The development lends itself to natural computational implementation for practical purposes.

Several positive aspects of the development in the context of Linear Elasticity have been pointed out. An unpalatable feature is that the transformation to an at most first-order differential system, due to the requirement that the pre-dual Lagrangian have no derivatives on primal variables appearing in it, increases the number of field degrees of freedom. The number of degrees of freedom in classical linear elastodynamics is 6 (we count velocity/momentum as a separate field as that is the optimal setting, especially when it comes to questions of heterogeneous materials, see e.g., (Willis 1980)). In contrast, the present formulation needs a minimum of 12 (6 for $\mu$ when Equation (2c) is posed in symmetrized form, and 3 each for displacement, and velocity). Thus, while providing conceptual benefits and exposing unexpected connections, the present approach is sub-optimal for standard problems of linear elasticity from a practical point of view. However, for non-standard or difficult problems, e.g., a convex variational formulation of 'Odd Elasticity' (Scheibner et al. 2020), problems with indefinite elastic moduli (when the dynamic Cauchy problem changes type and is ill-posed due to a lack of continuous dependence w.r.t. initial data), problems with non-unique solutions, or for heterogeneous materials with high contrast, the methodology can be potentially beneficial. 'Metamaterials' are included in this class of problems.

While the availability of a convex variational principle for a general PDE system is a definite advantage, the considerations related to degenerate ellipticity at the end of Section 3 make it clear that weak coercivity of the convex dual functional, i.e., roughly speaking, the growth of the functional to $+\infty$ as a relevant norm of its argument function tends to $\infty$, is not to be expected. Thus, one of the three criteria for use of the 'big hammer' of the Calculus of Variations (See Zeidler (1995) for the Generalized Weierstrass Theorem[4]) to prove existence of minimizers is not likely to be satisfied in many situations[5], when working within the present dual formulation. The formulation of Vorotnikov (2022), extending the pioneering work of Brenier (2018), is to be noted in this regard, but in any case, the gradient flow scheme proposed in (Vorotnikov and Acharya 2025) is perhaps the most pragmatic strategy to work with in efforts to obtain weak solutions of the original primal PDE through the dual scheme. That strategy is also immediately transferable to a computational algorithm for approximate solutions.

Finally, a speculation from a non-expert in the theory of homogenization: it seems like the natural transformation of the elastostatic(dynamic) problem of an arbitrary linear elastic composite to a space-time boundary value problem with constant-coefficient highest-order operator (38), along with a (formal) uniqueness guarantee in the usually encountered situations (major symmetry and positive definiteness for all phases), may be expected to provide some advantages for mathematical homogenization schemes for such materials.

[4] Strictly speaking, one needs the theorem for extended real-valued functionals in the present case.

[5] The affineness of the pre-dual $\widehat{S}$ in $D$ gives weak lower-semicontinuity, and one works in a closed, convex subset of a Hilbert space.

**Competing interests** The authors declare that they have no competing interests.

**Journal's Note** JTCAM remains neutral with regard to the content of the publication and institutional affiliations.